\documentclass[a4paper,twoside,reqno]{amsart}
\usepackage{amsthm,amsmath,amssymb,xcolor,soul,xspace,bm}
\usepackage[bookmarks=true,bookmarksopen=true]{hyperref}
\usepackage[ulem=normalem,draft]{changes}

\usepackage[ruled,algo2e]{algorithm2e}

\theoremstyle{plain}
\newtheorem{theorem}{Theorem}[section]
\newtheorem{lemma}[theorem]{Lemma}
\newtheorem{corollary}[theorem]{Corollary}

\theoremstyle{definition}
\newtheorem{definition}[theorem]{Definition}

\newtheorem{assumption}[theorem]{Assumption}

\theoremstyle{remark}
\newtheorem{remark}[theorem]{Remark}

\title[Semismooth Newton Method for Boundary Bilinear Control]{Semismooth Newton Method for Boundary Bilinear Control}

\thanks{The first and third authors were supported by MCIN/ AEI/10.13039/501100011033/ under research project PID2020-114837GB-I00.}

\thanks{The second author was supported by the Hellenic Foundation for Research and Innovation (H.F.R.I.) under the ``First Call for H.F.R.I. Research Projects to support Faculty members and Researchers and the procurement of high-cost research equipment grant" (Project Number: 3270)}

\author[E. Casas]{Eduardo Casas}
\address{Departamento de Matem\'{a}tica Aplicada y Ciencias de la Computaci\'{o}n, E.T.S.I. Industriales y de Telecomunicaci\'on, Universidad de Cantabria, 39005 Santander, Spain.}
\email{eduardo.casas@unican.es}

\author[K. Chrysafinos]{Konstantinos Chrysafinos}
\address{Department of Mathematics, School of Applied Mathematics and Physical Sciences, National Technical University of Athens, Zografou Campus, Athens 15780, Greece and IACM, FORTH, 20013 Heraklion, Crete, Greece.}
\email{chrysafinos@math.ntua.gr}

\author[M. Mateos]{Mariano Mateos}
\address{Departamento de Matem\'{a}ticas, Campus de Gij\'on, Universidad de Oviedo, 33203, Gij\'on, Spain.}
 \email{mmateos@uniovi.es}



\ifpdf
\hypersetup{
  pdftitle={Semismooth Newton Method for Boundary Bilinear Control},
  pdfauthor={E.~Casas,  K. Chrysafinos, and M.~Mateos}
}
\fi

\newcommand{\uad}{U_{\rm ad}}

\newcommand{\proj}{\mbox{\rm Proj}}
\newcommand{\dx}{\, \textup{d}x}
\newcommand{\dsigmax}{\, \textup{d}x}

\newcommand{\Pb}{\mbox{\rm (P)}\xspace}
\newcommand{\A}{\mathcal{A}}
\newcommand{\ci}{\chi_{{}_{\mathbb{I}_u}}\!}
\newcommand{\ca}{\chi_{{}_{\mathbb{A}_u}}\!}

\newcommand{\dimension}{d}
\newcommand{\conormal}{n}
\newcommand{\tichonov}{\nu}

\newcommand{\umin}{{\alpha}}
\newcommand{\umax}{{\beta}}

\begin{document}

\begin{abstract}
We study a control-constrained optimal control problem
governed by a semilinear elliptic equation. The control acts in a
bilinear way on the boundary, and can be interpreted as a heat
transfer coefficient. A detailed study of the state equation is
performed and differentiability properties of the control-to-state
mapping are shown. First and second order optimality conditions are
derived. Our main result is the proof of superlinear convergence of the
semismooth Newton method to local solutions satisfying no-gap
second order sufficient optimality conditions as well as a strict
complementarity condition. \end{abstract}

\maketitle

\section{Introduction}
\label{sec:introduction}
In this paper, we propose a semismooth Newton method to solve the following bilinear optimal control problem:
\[ \Pb \min_{u \in \uad} J(u) :=  \int_\Omega L(x,y_u(x)) \dx + \frac{\tichonov}{2} \int_\Gamma u^2(x) \dx, \]
where $y_u$ is the state associated with the control $u$ solution of

\begin{equation}
\left\{\begin{array}{l} Ay + a(x,y)= 0\ \  \mbox{in } \Omega,\vspace{2mm}\\  \partial_{\conormal_A} y +uy = g\ \ \mbox{on }\Gamma. \end{array}\right.
\label{E01}
\end{equation}
Here $\Omega \subset \mathbb{R}^\dimension$, $\dimension =2$ or $3$, is a bounded open connected set with a Lipschitz boundary $\Gamma$, $\tichonov >0$ and
\[ \uad = \{ u \in L^{2}(\Gamma) :  \umin \leq u(x) \leq \umax \text{ a.e. in }  \Gamma \}, \]
with $ 0 \leq  \umin < \umax < \infty.$
The remaining assumptions regarding the data of the control problem will be given in Sections \ref{S2} and \ref{S3}. Typical examples would include the tracking type functional $L(x,y)=\frac{1}{2}(y-y_d(x))^2$ for some target state $y_d$ and nonlinearities such as $a(x,y)=y^3$ or $a(x,y)=\exp(y)$.

Bilinear control plays an important role not only for the purposes of parameter identification, but also as ways of changing the intrinsic properties of the controlled system. Applications of bilinear control to very distinct fields such as nuclear and thermal control processes, ecologic and physiologic control or socioeconomic systems can be found in the early reference \cite{Mohler1973}, where they are investigated in the framework of ordinary differential equations. In the recent paper \cite{Winkler2020}, the author underlines the importance of bilinear boundary control of partial differential equations in several applications, providing references for them. The goal of that paper is not the analysis of an optimization algorithm, but the obtention of error estimates for the finite element approximation of \Pb, assuming that the state equation is linear.

Our main goal is to analyze the  convergence of the semismooth Newton method applied to \Pb. The novelty of this paper is twofold. First, the convergence analysis is carried out under the assumptions of no-gap second order optimality conditions and a strict complementarity condition, which are the usual ones to study numerical optimization algorithms in finite dimensional constrained optimization problems; see e.g. \cite{Nocedal-Wright1999}. This improves the previous results \cite{Amstutz-Laurain,Mannel-Rund2021,Pieper2015} for distributed controls and \cite{Reyes-Kunisch2005,Kroner-Kunisch-Vexler2011} for boundary controls, where conditions leading to local convexity were assumed. Second, as far as we know, there are no results in this direction for boundary bilinear controls. In \cite{Casas-Mateos2023} we considered a problem with distributed control acting as a source in the equation; in \cite{Casas-Chrysafinos-Mateos2023} we turned our attention to a bilinear control problem where the control appears as a reaction coefficient in the partial differential equation. In the paper at hand, the control appears as the Robin coefficient on the boundary condition and a new difficulty appears: the control-to-state mapping is not differentiable $L^2(\Gamma)$  if $\dimension=3$. In this paper, we focus on the aspects of the proofs that are essentially different from those in \cite{Casas-Chrysafinos-Mateos2023} and \cite{Casas-Mateos2023}, and refer to those papers when necessary.
\section{State equation} \label{S2}
Let us state the assumptions associated to the state equation.
\begin{assumption} \label{A2.1}
The operator $A$ is defined in $\Omega$ by
\[ Ay = - \sum_{i,j=1}^\dimension \partial_{x_j} [ a_{ij}(x) \partial_{x_i} y] +a_0y.\]
We suppose that $a_0,a_{ij} \in L^{\infty}(\Omega)$ for $1{\leq} i,j {\leq} \dimension$ with  $0\leq a_0\not\equiv  0$, and there exist $0<\tilde\lambda_A\leq\tilde\Lambda_A<\infty$ satisfying
\[\tilde\lambda_A | \xi |^2 \leq \sum_{i,j=1}^\dimension a_{ij}(x) \xi_i \xi_j \leq \tilde\Lambda_A | \xi |^2 \, \text{for a.a.} \ x {\in} \Omega \ \text{and} \ \forall \xi {\in} \mathbb R^\dimension.
\]

\end{assumption}
Notice that Assumption \ref{A2.1} implies the existence of $0<\lambda_A<\Lambda_A$ such that the bilinear form
\[\mathfrak{a}(y,z) = \int_\Omega \left(\sum_{i,j=1}^\dimension a_{ij}\partial_{x_i}y\partial_{x_j}z + a_0 y z\right)\dx \]
satisfies
\begin{align}
\mathfrak{a}(y,y)\geq &\lambda_A\|y\|_{H^1(\Omega)}^2&\forall y\in H^1(\Omega),\label{E02}\\
\mathfrak{a}(y,z)\leq &\Lambda_A\|y\|_{H^1(\Omega)}\|z\|_{H^1(\Omega)}&\forall y,z\in H^1(\Omega).\label{E03}
\end{align}

\begin{assumption}\label{A2.2}
We assume that $a:\Omega \times \mathbb{R} \longrightarrow \mathbb{R}$ is a Carath\'eodory function of class $C^2$ with respect to the second variable satisfying for a.a. $x\in\Omega$:
\begin{align*}
&\bullet a(\cdot,0) \in L^p(\Omega) \text{ for some } p>\dimension/2, \\
&\bullet  \frac{\partial a}{\partial y}(x,y) \geq 0\ \forall y\in\mathbb{R},\\
&\bullet \forall M {>} 0 \ \exists C_{a,M} \text{ s.t. } \sum_{j=1}^2 \Big| \frac{\partial^j a}{\partial y^j} (x,y) \Big| {\leq} C_{a,M}\, \forall |y|\leq M,\\
&\bullet \forall \varepsilon {>}0 \text{ and } \forall M{>}0 \, \exists \rho {>}0 \text{ s. t. }  \Big|    \frac{\partial^2 a}{\partial y^2}  (x,y_1)  -   \frac{\partial^2 a}{\partial y^2}  (x,y_2) \Big| {\leq} \varepsilon
\end{align*}
for all $ |y_1|, |y_2| \leq M$  with $|y_1-y_2| \leq \rho$.

All the above constants are independent of $x$.

We suppose that $g\in L^q(\Gamma)$ with $q> \dimension -1$ and, without loss of generality, that $q\leq d$.
\end{assumption}
To deal with the nonlinearity of the state equation, we observe that $q=2$ is not enough in dimension $d=3$. The proof of the differentiability of the relation control-to-state requires $q>2$. For linear state equations, $q=2$ is enough; see \cite{Winkler2020}.

For $\dimension = 2$ or $3$ it is
known that $H^{1/2}(\Gamma)\subset L^4(\Gamma)$ and there exists $C_\Gamma$ such that
\begin{equation} \|y\|_{L^{4}(\Gamma)} \leq C_{\Gamma} \|y\|_{H^1(\Omega)}, \quad\forall y \in H^1(\Omega). \label{E05}
\end{equation}

Throughout this paper the following notation will be used: we fix $s=2$ if $\dimension=2$ or $s=q$ if $\dimension =3$ and define the set
\begin{equation} \label{E06} \begin{array}{l}
\A_0 := \{ u \in L^{s}(\Gamma) : u \geq 0 \}.
\end{array}
\end{equation}
We denote $B_r(\bar u)=\{u\in L^s(\Gamma):\ \|u-\bar u\|_{L^s(\Gamma)}<r\}$.
\begin{theorem} \label{T2.3}
There exists $\mu>0$ such that for every $u \in \A_0$ equation \eqref{E01} has a unique solution $y_u \in Y:= H^1(\Omega) \cap C^{0,\mu}(\bar\Omega)$.
Furthermore, the following estimates hold:
\begin{align}
& \hskip-8pt \|y_u\|_{H^1(\Omega)} {\leq} C \left (  \|a(\cdot,0)\|_{L^p(\Omega)} + \|g\|_{L^q(\Gamma)} \right ), \label{E07} \\
& \hskip-8pt \|y_u\|_{L^\infty(\Omega)} \leq M_\infty(\|a(\cdot,0)\|_{L^p(\Omega)}+\|g\|_{L^q(\Gamma)}), \label{E08}\\
& \hskip-8pt \|y_u\|_{C^{0,\mu}(\bar\Omega)} \leq  C_{\mu,\infty}(\|a(\cdot,0)\|_{L^p(\Omega)}+\|u\|_{L^s(\Gamma)}+\|g\|_{L^q(\Gamma)}), \label{E09}
\end{align}
where $C$, $M_\infty$ and $C_{\mu,\infty}$ are independent of $u$.
\end{theorem}
\begin{proof}
We define the mapping
\begin{equation*} b: \Omega \times \mathbb R \longrightarrow \mathbb R, \   b(x,y) := a(x,y) - a(x,0).
\end{equation*}
Assumption \ref{A2.2} implies that $b(x,0)=0$ and $\frac{\partial b}{\partial y} (x,y) \geq 0$.
Equation \eqref{E01} can be written in the variational form
\begin{align}
&\mathfrak{a}(y,z) + \int_\Omega b(x,y)z \dx +\int_\Gamma u y z \dsigmax\nonumber \\
&=
\int_\Omega -a(x,0) z \dx +\int_\Gamma g z \dsigmax\quad \forall z\in H^1(\Omega).
 \label{E10}
 \end{align}
Using  \eqref{E03}, Cauchy's inequality, \eqref{E05}, \eqref{E02} and the nonnegativity of $u$ imposed in \eqref{E06}, we infer that
\begin{align}
\mathfrak{a}(y,z)+\int_\Gamma u y z \dsigmax \leq& \Lambda_u\|y\|_{H^1(\Omega)}\|z\|_{H^1(\Omega)}, \label{E11}\\
\mathfrak{a}(y,y)+\int_\Gamma u y^2 \dsigmax \geq &\lambda_A\|y\|_{H^1(\Omega)}^{2}, \label{E12}
\end{align}
where $\Lambda_u {=} \Lambda_A + \|u\|_{L^{2}(\Gamma)} C_\Gamma^2$. The proof of existence and uniqueness of a solution in $H^1(\Omega)\cap L^\infty(\Omega)$ of \eqref{E10} as well as estimates \eqref{E07} and  \eqref{E08} follow as in \cite[Theorem 3.1]{Casas93}. The $L^\infty(\Omega)$ estimate is obtained following the approach of \cite[Theorem 4.1]{Stampacchia65} and using that $u \ge 0$ and $b(x,s)s \ge 0$ $\forall s \in \mathbb{R}$.

To prove \eqref{E09}  we write \eqref{E01} in the form
\begin{equation*}
\left\{\begin{array}{l} Ay = - a(x,y) \ \  \mbox{in } \Omega,\vspace{2mm}\\  \partial_{\conormal_A} y  = -uy + g\ \ \mbox{on }\Gamma. \end{array}\right.
\end{equation*}
From Assumption \ref{A2.2} and the mean value theorem we infer
\[ |a(x,y) | \leq | a(x,0) | + C_{a,K}K, \]
where $K = \|y\|_{L^\infty(\Omega)}$. In addition, we have $\|-u y\|_{L^s(\Gamma)} \leq  K \|u\|_{L^s(\Gamma)}$. Then, from \cite[Proposition 3.6]{Nittka11} we infer that $y$ belongs to  $C^{0,\mu}(\bar \Omega)$ and satisfies \eqref{E09} for some $\mu \in (0,1]$.
\end{proof}
Next we consider the differentiability of the mapping $u \to y_u$.
\begin{theorem} \label{T2.4}
There exists an open set $\A$ in $L^{s}(\Gamma)$ such that $\A_0 \subset \A$ and equation \eqref{E01} has a unique solution $y_u \in Y$ $\forall u \in \A$.
Further, the mapping $G : \A \longrightarrow Y$ defined by $G(u):=y_u$ is of class $C^2$ and $\forall u \in \A$ and $\forall v, v_1, v_2 \in L^{s}(\Gamma)$ the functions
$z =G'(u)v$ and $w=G''(u)(v_1,v_2)$ are the unique solutions of the equations:
\begin{align}
& \left\{\begin{array}{l} Az + \displaystyle \frac{\partial a}{\partial y}(x,y_u)z = 0\ \  \mbox{in } \Omega,\vspace{2mm}\\  \partial_{\conormal_A} z + uz  = -v y_u\ \ \mbox{on }\Gamma, \end{array}\right.
\label{E13} \\
& \left\{\begin{array}{l} Aw + \displaystyle \frac{\partial a}{\partial y}(x,y_u)w+ \displaystyle \frac{\partial^2 a}{\partial y^2}(x,y_u)z_{u,v_1}z_{u,v_2} = 0\ \  \mbox{in } \Omega,\vspace{2mm}\\  \partial_{\conormal_A} w + uw = - v_1z_{u,v_2} - v_2z_{u,v_1} \ \ \mbox{on }\Gamma, \end{array}\right.
\label{E14}
\end{align}
where $z_{u,v_i} = G'(u)v_i$, $i=1,2$.
\end{theorem}
\begin{proof}
We consider the space
\[Y_A := \{ y \in Y : Ay \in L^p(\Omega) , \  \partial_{\conormal_A} y \in L^q(\Gamma) \}  \]
endowed with the graph norm. We note that $Y_A$ is a Banach space.
We also define the mapping ${\mathcal F}: L^{s}(\Gamma) \times Y_A \longrightarrow L^p(\Omega) \times L^q(\Gamma)$ by
\[ {\mathcal F}(u,y) := ( Ay + a(\cdot,y), \partial_{\conormal_A} y +uy- g ). \]
Since $q\leq s$, $\mathcal F$ is well defined and of class $C^2$ due to Assumption \ref{A2.2}. For every $(u, y)\in\A_0\times Y_A$ the derivative $\frac{\partial {\mathcal F}}{\partial y}( u, y) : Y_A \longrightarrow L^p(\Omega) \times L^q(\Gamma)$, given by
\[ \frac{\partial {\mathcal F}}{\partial y}( u, y)z = \left  ( Az + \displaystyle \frac{\partial a}{\partial y}(\cdot, y)z,    \partial_{\conormal_A} z +  u z\right ) \ \forall z \in Y_A,\]
is linear and continuous.
The open mapping theorem implies that  $\frac{\partial {\mathcal F}}{\partial y}( u, y)$ is an isomorphism if and only if the equation,
 \begin{equation*}
 \left\{\begin{array}{l} Az + \displaystyle \frac{\partial a}{\partial y}(x, y)z = f\ \  \mbox{in } \Omega,\vspace{2mm}\\  \partial_{\conormal_A} z +  u z = h\ \ \mbox{on }\Gamma, \end{array}\right.
\end{equation*}
has unique solution $z \in Y_A$ for all $(f,h) \in L^p(\Omega) \times L^q(\Gamma).$ This fact follows from Theorem \ref{T2.3}.
Then, given $\bar u\in\A_0$ with $\bar y=y_{\bar u}$, since $\mathcal F (\bar u,\bar y) = 0$, the implicit function theorem implies the existence of $ \varepsilon_{\bar u} >0$ and $\varepsilon_{\bar y} >0$ such that $\forall u \in B_{\varepsilon_{\bar u}}(\bar u)\subset L^s(\Gamma)$ the equation ${\mathcal F}(u,y) = 0$ has a unique solution $y_u$ in the open ball $B_{\varepsilon_{\bar y}}(\bar y) \subset Y_A\subset Y.$  Moreover, the mapping $u \in  B_{\varepsilon_{\bar u}}(\bar u) \to y_u \in B_{\varepsilon_{\bar y}}(\bar y)$ is of class $C^2$. Without loss of generality, we assume $\varepsilon_{\bar u} < \frac{1}{2}\lambda_A/(|\Gamma|^{\frac{s-2}{s}}C_\Gamma^2)$, where $C_\Gamma$ is introduced in \eqref{E05}.
Actually, for every $u \in B_{\varepsilon_{\bar u}}(\bar u)$ the equation ${\mathcal F}(u,y) = 0$ has unique solution $y \in Y_A$. Indeed, let $y_1,y_2$ denote two solutions of
${\mathcal F}(u,y) = 0.$ We set $y=y_2-y_1$, subtract the corresponding equations, and apply the mean value theorem to deduce that $y$ satisfies
 \begin{equation} \label{E15}
 \left\{\begin{array}{l} Ay + \displaystyle \frac{\partial a}{\partial y}(x,y_1+\theta_x y) y = 0\ \  \mbox{in } \Omega,\vspace{2mm}\\  \partial_{\conormal_A} y + u y= 0 \ \ \mbox{on }\Gamma, \end{array}\right.
\end{equation}
where $\theta_x:\Omega\to[0,1]$ is a measurable function.
Adding and subtracting appropriate terms on the boundary, equation \eqref{E15} can be written as
 \begin{equation} \label{E16}
 \left\{\begin{array}{l} Ay + \displaystyle \frac{\partial a}{\partial y}(x,y_1+\theta_x y) y = 0\ \  \mbox{in } \Omega,\vspace{2mm}\\  \partial_{\conormal_A} y +\bar u y =- (u-\bar u)y \ \ \mbox{on }\Gamma. \end{array}\right.
\end{equation}
Testing the variational form of \eqref{E16} with $y$ we get
\begin{align*}
 \lambda_A \|y\|^2_{H^1(\Omega)}
&  \leq  \varepsilon_{\bar u} |\Gamma|^{\frac{s-2}{s}}C^2_{\Gamma}  \|y\|^2_{H^1(\Omega)}.
\end{align*}
Since $\varepsilon_{\bar u} < \frac{1}{2}\lambda_A/(|\Gamma|^{\frac{s-2}{s}}C^2_{\Gamma})$, $y=0$ holds. Defining in $L^{s}(\Gamma)$ the open set $\A = \cup_{\bar u \in \A_0} B_{\varepsilon_{\bar u}} (\bar u)$ and $G: \A \longrightarrow Y$ such that $G(u)=y_u,$ we have that $G$ is of class of $C^2$. Finally, equations \eqref{E13} and \eqref{E14} are obtained differentiating with respect to $u$ the identity ${\mathcal F}(u,G(u)) = 0$.
\end{proof}
\begin{remark}\label{R2.5}
Theorems \ref{T2.3} and \ref{T2.4} are valid if we use the operator $A^*$ instead of $A$, where
$ A^*\varphi = - \sum_{i,j=1}^\dimension \partial_{x_j} [ a_{ji}(x) \partial_{x_i} \varphi] +a_0\varphi$.
Therefore, for every $\bar u\in\A_0$ we obtain the existence of $\varepsilon_{\bar u}^*>0$ such that, for every $(f,h)\in L^p(\Omega)\times L^q(\Gamma)$ and $u\in B_{\varepsilon_{\bar u}^*}(\bar u)$, the equation
\[
 \left\{\begin{array}{l} A^*\varphi + \displaystyle \frac{\partial a}{\partial y}(x,y_u)\varphi  = f \ \  \mbox{in } \Omega,\vspace{2mm}\\  \partial_{\conormal_{A^*}} \varphi + u\varphi = h \ \ \mbox{on }\Gamma, \end{array}\right.
\]
has a unique solution $\varphi\in Y$. Without loss of generality, we can assume that $\varepsilon_{\bar u}\leq\varepsilon_{\bar u}^*$, so the equation is uniquely solvable in $Y$ for all $u\in\A$.
\end{remark}

\section{Analysis of the optimal control problem} \label{S3}
In this section we proceed to the analysis of the optimal control problem. To this end we make the following hypotheses on $J$.

\begin{assumption} \label{A3.1} The function $L: \Omega \times \mathbb R \longrightarrow \mathbb R$ is Carath\'eodory and of class of $C^2$ with respect to the second variable.
Further the following properties hold for a.a. $x \in \Omega$:
\begin{align*}
& \bullet L(\cdot, 0) \in L^1(\Omega), \\
& \bullet \forall M >0, \  \exists L_M \in L^p(\Omega) \text{ such that } \Big|  \frac{\partial L}{\partial y}(x,y) \Big| \leq L_M(x), \\
& \bullet \forall M >0, \ \exists C_{L,M} \in \mathbb R   \text{ such that } \Big|  \frac{\partial^2 L}{\partial y^2}(x,y) \Big| \leq C_{L,M}, \\
& \bullet \forall \varepsilon >0  \text{ and } \forall M >0 \ \exists \rho >0 \text{ such that }
\Big|  \frac{\partial^2 L}{\partial y^2}(x,y_1) - \frac{\partial^2 L}{\partial y^2}(x,y_2)  \Big| \leq \varepsilon
\end{align*}
for all $|y|, |y_1|, |y_2| \leq M$  with $|y_1-y_2| \leq \rho$.
All the above constants are independent of $x$.
\end{assumption}
The following theorem states the differentiability properties of the minimizing functional.

\begin{theorem} \label{T3.2} The functional $J: \A \longrightarrow \mathbb R$ is of class $C^2$ and its derivatives are given by the expressions:
\begin{align}
 J'(u)v = &\int_\Gamma (\tichonov u - y_u \varphi_u)v  \dsigmax, \label{E17} \\
 J''(u)(v_1,v_2) =
 \nonumber
  & \int_\Omega \Big[ \frac{\partial^2 L}{\partial y^2}(x,y_u) - \varphi_u  \frac{\partial^2 a}{\partial y^2}(x,y_u) \Big] z_{u,v_1} z_{u,v_2}\dx
  \\ \label{E18}
&
- \int_\Gamma \Big[ v_1z_{u,v_2} + v_2 z_{u,v_1} \Big] \varphi_u \dsigmax + \tichonov \int_\Gamma v_1 v_2 \dsigmax,
\end{align}
for all  $u \in \A$ and all $v, v_1, v_2 \in L^{s}(\Gamma)$,
where $z_{u,v_i} = G'(u)v_i$, $i=1,2$ and $\varphi_u \in Y$ is the adjoint state, the unique solution of the equation
\begin{equation}
 \left\{\begin{array}{l} A^*\varphi + \displaystyle \frac{\partial a}{\partial y}(x,y_u)\varphi  = \frac{\partial L}{\partial y}(x,y_u) \ \  \mbox{in } \Omega,\vspace{2mm}\\  \partial_{\conormal_{A^*}} \varphi + u\varphi = 0\ \ \mbox{on }\Gamma. \end{array}\right.
\label{E19}
\end{equation}
\end{theorem}
\begin{proof}
Existence, uniqueness and regularity of $\varphi_u$ follow from Remark \ref{R2.5}, Assumption \ref{A3.1}, and Theorem \ref{T2.4}.
The proofs of \eqref{E17} and \eqref{E18} are standard and can be established working identically to \cite[Theorem 3.4]{Casas-Chrysafinos-Mateos2023}.
\end{proof}

According to Theorem \ref{T3.2} the mapping $\Phi : \A \longrightarrow Y$ given by $\Phi (u):= \varphi_u$ is well defined. Let us prove that it is $C^1$.
\begin{theorem} \label{T3.3}
The mapping $\Phi$ is of class $C^1$ and for all $u\in \A$ and $v\in L^{s}(\Gamma)$ the function $\eta_{u,v} = \Phi'(u)v$ is the unique solution of
\begin{equation}
 \left\{\begin{array}{l}
 A^*\eta {+} \displaystyle \frac{\partial a}{\partial y}(x,y_u)\eta  = \Big[\frac{\partial^2 L}{\partial y^2}(x,y_u) - \varphi_u  \frac{\partial^2 a}{\partial y^2}(x,y_u)\Big] z_{u,v} \,  \mbox{in}\, \Omega,\vspace{2mm}\\
 \!\!\!\! \partial_{\conormal_{A^*}} \eta {+} u\eta {=} - v\varphi_u \ \ \mbox{on }\Gamma, \end{array}\right.
\label{E21}
\end{equation}
where $z_{u,v}= G'(u)v.$
\end{theorem}
\begin{proof} Using Assumption \ref{A3.1} and the fact that $y_u,\varphi_u,z_{u,v} \in L^{\infty}(\Omega)$ we obtain that the right hand side of \eqref{E21} belongs to $L^{p}(\Omega)\times L^s(\Gamma)$. Existence, uniqueness, and regularity of $\eta_{u,v}$ follow from Remark \ref{R2.5}.
To establish the differentiability of $\Phi$ we define
\begin{align*}
Y_{A^*}  = \{ \varphi \in Y: A^* \varphi \in L^p(\Omega)  \text{ and } \partial_{\conormal_{A^*}} \varphi \in L^q(\Gamma) \}
\end{align*} and
$\mathcal G : \A \times Y_{A^*} \longrightarrow L^p(\Omega) \times L^q(\Gamma)$ by
\[ {\mathcal G}(u,\varphi ) := \Big(A^* \varphi + \frac{\partial a}{\partial y}(\cdot,y_u) \varphi - \frac{\partial L}{\partial y}(\cdot,y_u),\partial_{\conormal_{A^*}} \varphi + u \varphi \Big).  \]
From assumptions \ref{A2.2} and \ref{A3.1} we deduce that $\mathcal G$ is of class $C^1$.  Moreover, $\frac{\partial \mathcal G}{\partial \varphi}(u,\varphi): Y_{A^*} \longrightarrow L^p(\Omega) \times L^q(\Gamma)$ is a linear and continuous mapping, and $\forall \eta \in Y_{A^*}$ we have  that
\[
\frac{\partial \mathcal G}{\partial \varphi}(u,\varphi) \eta = \left ( A^* \eta + \frac{\partial a}{\partial y}(\cdot,y_u) \eta, \partial_{\conormal_{A^*}} \eta + u \eta \right ).
\]
Using again Remark \ref{R2.5} we get that
\begin{equation*}
 \left\{\begin{array}{l} A^*\eta + \displaystyle \frac{\partial a}{\partial y}(x,y_u)\eta = f\ \  \mbox{in } \Omega,\vspace{2mm}\\  \partial_{\conormal_{A^*}}\eta + u\eta = h\ \ \mbox{on }\Gamma, \end{array}\right.
\end{equation*}
has a unique solution in $Y_{A^*}$ for all $(f,h) \in L^p(\Omega) \times L^q(\Gamma)$. Hence, $\frac{\partial \mathcal G}{\partial \varphi}(u,\varphi) : Y_{A^*} \longrightarrow L^p(\Omega) \times L^q(\Gamma)$ is an isomorphism. Then, applying the implicit function theorem and differentiating the identity $\mathcal{G}(u,\Phi(u)) = 0$ the result follows.
\end{proof}

Combining \eqref{E21} with \eqref{E18} we deduce the following alternative representation formula for $J''(u)$.
\begin{corollary} \label{C3.4} For every $v_1, v_2 \in L^{s}(\Gamma)$ and all $u \in \A$, the following identities hold
\begin{align}
J''(u)(v_1,v_2)  = &\int_{\Gamma}\Big[ \tichonov v_1 - (\varphi_u z_{u,v_1}  + y_u \eta_{u,v_1}) \Big] v_2 \dsigmax
=  \int_{\Gamma}\Big[ \tichonov v_2 - (\varphi_u z_{u,v_2}  + y_u \eta_{u,v_2}) \Big] v_1 \dsigmax.
\label{EX22a}
\end{align}
\end{corollary}

\begin{remark}\label{R3.5}In dimension $\dimension=3$, we can also extend $J'(u)$ and $J''(u)$ respectively to continuous linear and bilinear forms in $L^2(\Gamma)$ and $L^2(\Gamma)^2$ by the same expressions given above. Indeed, we notice that for all $v\in L^2(\Gamma)$, the Lax-Milgram Theorem implies that equations \eqref{E13} and \eqref{E21} have a unique solution in $H^1(\Omega) \subset L^2(\Omega)$.
\end{remark}

\begin{theorem} \label{T3.6} Problem $\Pb$ has at least one solution. Moreover, if $\bar u \in \uad$ is a local minimizer of $\Pb$ then there exist $\bar y, \bar \varphi \in Y$ such that
\begin{align}
& \left\{\begin{array}{l} A\bar y + a(x,\bar y) = 0\ \  \mbox{in } \Omega,\vspace{2mm}\\  \partial_{\conormal_A} \bar y + \bar u \bar y = g\ \ \mbox{on }\Gamma, \end{array}\right.
\label{E22} \\
& \left\{\begin{array}{l} A^* \bar \varphi + \displaystyle \frac{\partial a}{\partial y}(x,\bar y)\bar \varphi = \frac{\partial L}{\partial y}(x,\bar y) \ \  \mbox{in } \Omega,\vspace{2mm}\\  \partial_{\conormal_{A^*}} \bar  \varphi + \bar u \bar \varphi = 0\ \ \mbox{on }\Gamma, \end{array}\right. \label{E23} \\
& \bar u(x) = \proj_{[\umin, \umax]} \left ( \frac{1}{\tichonov} \bar y(x) \bar \varphi(x) \right )\qquad \forall x\in\Gamma. \label{E24}
\end{align}
Moreover, the regularity $\bar u \in C^{0,\mu}(\Gamma)$ holds.
\end{theorem}
\begin{proof}
Existence of optimal solutions follows by the direct method of the calculus of variations. The only delicate point is to show that for every sequence $\{u_k\}_{k=1}^\infty\subset \uad$ such that $u_k\rightharpoonup u$ weakly in $L^s(\Gamma)$, the sequence $y_{u_k}\to y_u$ converges strongly in $C(\bar\Omega)$. From Theorem \ref{T2.3} we have that $\{y_{u_k}\}_k$ is bounded in $Y$. Hence, there exists $y\in Y$ such that $y_{u_k}\rightharpoonup y$ weakly in $H^1(\Omega)$. The compactness of the embedding $C^{0,\mu}(\bar\Omega)\subset C(\bar\Omega)$ implies the strong convergence in $C(\bar\Omega)$  and, consequently, $u_ky_{u_k}\rightharpoonup uy$ weakly in $L^s(\Gamma)$. Therefore, we can take limits in the equation satisfied by $y_{u_k}$ to deduce that $y = y_u$.

First order optimality conditions are an immediate consequence of \eqref{E17} and the convexity of $\uad$. The H\"older continuity of $\bar u$ is a
consequence of \eqref{E24}, the same regularity for $\bar y$ and $\bar\varphi$, and the Lipschitz property of the projection
$\proj_{[\umin, \umax]}(t) = \max\{\umin,\min\{\umax,t\}\}$.
\end{proof}
In this paper a local minimizer is intended in the $L^2(\Gamma)$ sense.
From now on $(\bar u,\bar y,\bar \varphi) \in \uad \times Y^2$ will denote a triplet satisfying \eqref{E22}-\eqref{E24}.
Associated with this triplet we define the cone of critical directions
\begin{equation*}  \hskip-12pt
\begin{array}{l} C_{\bar u} {=}
\{v {\in} L^{2}(\Gamma) :\, v(x) {=} 0 \text{ if } \tichonov \bar u(x)  -  \bar y(x) \bar \varphi(x) {\ne} 0
\text{ a.e. in $\Gamma$ and (\ref{E25})  holds} \},
\end{array}
\end{equation*}
\begin{equation} \label{E25}
v(x)  \left\{ \begin{array}{ll} \geq 0 & \text{ if } \bar u(x) = \umin, \\ \leq 0 & \text{ if } \bar u(x) = \umax. \end{array} \right.
\end{equation}

We proceed now to the second order optimality conditions. The proof of the following theorem is standard; see, e.g. \cite[Theorem 2.3]{Casas-Troltzsch12}.
\begin{theorem} \label{T3.7}
If $\bar u$ is a local minimizer of $\Pb$, then $J''(\bar u)v^2 \geq 0 \ \forall v \in C_{\bar u}$ holds. Conversely, if $\bar u \in \uad$ satisfies the first order optimality conditions \eqref{E22}--\eqref{E24} and $J''(\bar u)v^2 > 0 \ \forall v \in C_{\bar u} \setminus\{0\}$, then there exist $\varepsilon >0$ and $\delta >0$ such that
\begin{equation*}
J(\bar u) + \frac{\delta}{2} \|u - \bar u\|^2_{L^2(\Gamma)} \leq J(u)\, \forall u \in \uad\, \text{with}\, \|u-\bar u\|_{L^2(\Gamma)} \leq \varepsilon.
\end{equation*}
\end{theorem}

\smallskip
\begin{definition} \label{D3.8}
Let us define
\[ \Sigma_{\bar u} = \{ x \in \Gamma : \bar u(x) \in \{ \umin, \umax \} \text{ and } \tichonov \bar u(x) - \bar y(x) \bar \varphi (x) = 0 \}.\]
We say that the strict complementarity condition is satisfied at $\bar u$ if $ | \Sigma_{\bar u} | = 0,$ where $ | \cdot | $ stands for the $(\dimension-1)$ dimensional Lebesgue measure on $\Gamma$.
\end{definition}

For every $\tau \geq 0$, we define the subspace
\begin{equation*}
 T^{\tau}_{\bar u} {=} \{ v {\in} L^2(\Gamma) : v(x) {=} 0 \text{ if } | \tichonov \bar u(x)  -  \bar y (x) \bar\varphi (x) | > \tau\}.
\end{equation*}
\begin{theorem} \label{T3.9}
Assume that $\bar u$ satisfies the strict complementarity condition. Then, the following properties hold: \\
1- $T^0_{\bar u} = C_{\bar u}.$ \\
2- If $\bar u$ satisfies the second order optimality condition $J''(\bar u)v^2 >0$ $\forall v \in C_{\bar u} \setminus\{0\}$, then $\exists \tau >0$ and $\kappa >0$ such that
\begin{equation} \label{E26}
J''(\bar u) v^2 \geq \kappa \|v\|^2_{L^2(\Gamma)} \ \forall v \in T^\tau_{\bar u}.
\end{equation}
\end{theorem}
For the proof the reader is referred to \cite[Theorem 3.10]{Casas-Chrysafinos-Mateos2023}.

\section{Convergence of the semismooth Newton method} \label{S4}

We define $F{:} \A {\longrightarrow} L^{s}(\Gamma)$  by
$ F(u) {=} u - \proj_{[\umin,\umax]} \left (\frac{1}{\tichonov} y_u \varphi_u \right )$.
From theorems \ref{T2.4} and \ref{T3.2} we deduce that $F$ is well defined.
Due to Theorem \ref{T3.6}, any local minimizer of $\Pb$ is a solution of $F(u) = 0$. If a local minimizer $\bar u$ satisfies
$J''(\bar u)v^2 >0$ $\forall v \in C_{\bar u} \setminus\{0\}$, then there exists $\delta > 0$ such that it is the unique stationary point of $J$ in $B_\delta(\bar u)\cap\uad$; see \cite[Corollary 2.6]{Casas-Troltzsch12}. We are going to apply the semismooth Newton method sketched in Algorithm \ref{Alg1} to solve this equation.
\LinesNumbered
\begin{algorithm2e}[h!]
\caption{Semismooth Newton method.}\label{Alg1}
\DontPrintSemicolon
Initialize Choose $u_0\in \A$. Set $j=0$.\;
\For{$j \geq 0$}{
Choose $M_j\in\partial F(u_j)$ and solve $M_jv_j=-F(u_j)$\label{line3}.\;
Set $u_{j+1}=u_j+v_j$ and $j=j+1$.\;
}
\end{algorithm2e}
Here $\partial F(u)$ is a set valued mapping such that $F$ is $\partial F$-semismooth in the sense stated in \cite[Chapter 3]{Ulbrich11}. Local superlinear convergence follows from
the semismoothness of $F$ and the uniform boundedness of the norms of the inverses of the operators $M_j$.
In order to define $\partial F(u) \ \forall u \in \A$ we introduce some additional functions.
\begin{align*}
& S : \A  \longrightarrow L^{s}(\Gamma), \quad S(u)= \frac{1}{\tichonov} G(u) \Phi (u), \\
& \psi : \mathbb R  \longrightarrow \mathbb R, \quad  \psi(t) = \proj_{[\umin,\umax]} (t), \\
& \Psi :  \A \longrightarrow  L^{s}(\Gamma), \quad  \Psi(u)(x) = \psi(S(u)(x)).
\end{align*}
For every $u \in \A$ we define
\begin{align*} \partial \Psi (u) = & \big \{ N \in {\mathcal L}(L^{s}(\Gamma),L^{s}(\Gamma))\!:\! Nv = h S'(u)v \ \forall v \in L^{s}(\Gamma)\\
& \text{and for some measurable  function}\\
& h: \Omega \longrightarrow \mathbb R \text{ such that }  h(x) \in \partial \psi (S(u)(x))\big \}.
\end{align*}
We observe that $\psi$ is a Lipschitz function and by $\partial \psi (t)$ we denote the subdifferential in Clarke's sense; see \cite[Chapter 2]{Clarke83}. Note that
\[ \partial \psi(t) = \left\{ \begin{array}{cl} \{1\} & \text{ if } t \in (\umin,\umax), \\ \{0\} & \text{ if } t \not\in [\umin, \umax], \\  {[0,1]} & \text{ if } t \in \{\umin,\umax \}.
 \end{array}\right. \]
According to \cite[Prop. 2.26]{Ulbrich11}, $\psi$ is $1$-order $\partial \psi$-semismooth.

\begin{theorem} \label{T4.1}
$\Psi$ is $\partial \Psi$-semismooth in $\A$.
\end{theorem}

\begin{proof}
  Since $\Psi$ is a superposition operator of $\psi$ and $S$, we will apply \cite[Theorem 3.49]{Ulbrich11} to deduce that $\partial \Psi$-semismooth in $\A$.
To this end it is enough to prove that $ S : \A  \longrightarrow L^s(\Gamma)$ is $C^1$ and that $ S : \A  \longrightarrow L^r(\Omega)$ is locally Lipschitz for some $r>s$. The first condition is an immediate consequence of Theorems \ref{T2.4} and \ref{T3.3}. Indeed, since $S(u) = \frac{1}{\tichonov} G(u) \Phi(u)$ we have that
\[S'(u)= \frac{1}{\tichonov} [ G'(u)v \Phi(u) + G(u) \Phi'(u)v ] = \frac{1}{\tichonov} [ z_{u,v}\varphi_u + y_u \eta_{u,v} ]. \]
The Lipschitz condition is an immediate consequence of Lemma \ref{L4.2}.
\end{proof}

\begin{lemma}\label{L4.2}
For all $\bar u\in\A_0$, there exists $L_S>0$ such that
  \[\|S(u_1)-S(u_2)\|_{C(\Gamma)}\leq L_S \|u_1-u_2\|_{L^s(\Gamma)}\ \forall u_1,u_2\in B_{\varepsilon_{\bar u}}(\bar u)\]
  where $\varepsilon_{\bar u}$ is the one introduced in Theorem \ref{T2.4}.
\end{lemma}

\begin{proof}
  This is a consequence of Lemmas \ref{LA.2}, \ref{LA.3}, \ref{LA.7}, and \ref{LA.8}:
  \begin{align*}
    \|S&(u_1)-S(u_2)\|_{C(\Gamma)} \leq \|y_{u_1}\varphi_{u_1}-y_{u_2}\varphi_{u_2}\|_{C(\bar\Omega)} \\
    \leq &     \|y_{u_1}(\varphi_{u_1}-\varphi_{u_2})\|_{C(\bar\Omega)}+
        \|(y_{u_1}-y_{u_2})\varphi_{u_2}\|_{C(\bar\Omega)}
        \leq (K_\infty L_{\Phi,u}+K_\infty^* L_{G,u})\|u_1-u_2\|_{L^s(\Gamma)}.
  \end{align*}
\end{proof}

\begin{corollary}
The function $F:\A \longrightarrow L^{s}(\Gamma)$ is $\partial F$-semismooth in $\A$, where
\[
\partial F(u) = \{M = I - N : N \in \partial\Psi(u)\}
\]
and $I$ denotes the identity in $L^s(\Gamma)$.
\label{C4.3}
\end{corollary}

We select the operators $M_u:L^{s}(\Gamma) \longrightarrow L^{s}(\Gamma)$ for every $u \in \A$ as follows. First, we define the function $\lambda:\mathbb{R} \longrightarrow \mathbb{R}$ by
\[
\lambda(t) = \Big\{\begin{array}{cl} 1&\text{if } t \in (\umin,\umax),\\0&\text{otherwise}.\end{array}
\]
It is obvious that $\lambda(t) \in \partial\psi(t)$ for every $t \in \mathbb{R}$. We define $M_u:L^{s}(\Gamma) \longrightarrow L^{s}(\Gamma)$ by $M_uv = v - h_u\cdot S'(u)v$, where $h_u(x) = \lambda(S(u)(x)) = \lambda\big(\frac{1}{\tichonov}y_u(x)\varphi_u(x)\big)$. It is immediate that $M_u \in \partial F(u)$. For this selection we have the following result.

\begin{theorem}
Let $(\bar u,\bar y,\bar\varphi) \in \uad \times Y^2$ satisfy the first order optimality conditions \eqref{E22}--\eqref{E24}, the strict complementarity condition $|\Sigma_{\bar u}| = 0$, and the second order sufficient optimality condition $J''(\bar u)v^2 > 0$ for every $v \in C_{\bar u}\setminus\{0\}$. Then, there exist $\delta > 0$ and $C > 0$ such that for every $u \in B_\delta(\bar u) \subset \A$ the linear operator $M_u:L^{s}(\Gamma) \longrightarrow L^{s}(\Gamma)$ is an isomorphism and $\|M_u^{-1}\| \le C$ holds.
\label{T4.4}
\end{theorem}

\begin{proof}
For any $u\in\A$, we define
\[
\mathbb{A}_u = \{x \in \Gamma : \frac{1}{\tichonov}y_u(x)\varphi_u(x) \not\in (\umin,\umax)\} \ \text{ and }\ \mathbb{I}_u = \Gamma \setminus \mathbb{A}_u.
\]
Thus, the identity $M_uv = v - \frac{1}{\tichonov} [z_{u,v}\varphi_u + y_u \eta_{u,v}]\ci$ holds. Here $\chi_{\mathbb{S}}$ stands for the characteristic function of a set $\mathbb{S}$. $M_u$ being obviously continuous. Then, as a consequence of the open mapping theorem, it is enough to prove that the equation $M_u v = w$ has a unique solution $v\in L^s(\Gamma)$ for every $w\in L^{s}(\Gamma)$ to infer that $M_u$ is an isomorphism. Clearly, $v=w$ in $\mathbb{A}_u$, and hence, denoting $b = w + \frac{1}{\tichonov} [z_{u,\ca w}\varphi_u + y_u\eta_{u,\ca w}]\in L^s(\Gamma)$, to compute $v$ we have to solve
\begin{equation}
\ci v - \frac{1}{\tichonov} [z_{u,\ci v}\varphi_u + y_u\eta_{u,\ci v}] = b\text{ in }\mathbb{I}_u.
\label{E28}
\end{equation}
Using \eqref{EX22a}, it is obvious that this equation is the optimality condition of the unconstrained quadratic optimization problem
\[(Q) \min_{v\in L^2(\mathbb{I}_u)} \mathbb{J}(v)= \frac{1}{2\tichonov}J''(u)(\ci v)^2 - \int_{\mathbb{I}_u} b v \dsigmax.\]
Here and elsewhere, for every measurable set $\Sigma\subset\Gamma$ and $v\in L^1(\Sigma)$, $\chi_{_\Sigma}v$ denotes the extension by $0$ to $\Gamma\setminus\Sigma$.
The continuity of $J''$ established in Lemma \ref{LA.11} and \eqref{E26} imply the existence of  $\delta_0>0$ such that
\begin{equation}\label{E27}
J''(u)v^2 \ge \frac{\kappa}{2}\|v\|^2_{L^2(\Gamma)}\quad \forall v \in T^\tau_{\bar u}\text{ and } \forall u \in B_{\delta_0}(\bar u).
\end{equation}
Setting $\delta = \min\{\delta_0,\varepsilon_{\bar u}, \frac{\tau}{\tichonov L_S}\}$, where $\varepsilon_{\bar u}$  and $L_S$ are introduced in Theorem \ref{T2.4} and Lemma \ref{L4.2} respectively, we have that $L^2(\mathbb{I}_u)\subset T^\tau_{\bar u}$ for all  $u\in B_\delta(\bar u)$. To check this, we have to prove that
  \[\mathbb{I}_u\subset \{x\in\Gamma: |\tichonov \bar u(x)- \bar y(x)\bar\varphi(x)|\leq \tau\},\]
  or, equivalently, that
  \[ \{x\in\Gamma: |\tichonov \bar u(x)- \bar y(x)\bar\varphi(x)| > \tau\}\subset \mathbb{A}_u.\]
  If $\tichonov\bar u(x)- \bar y(x)\bar\varphi(x) > \tau$, then the first order optimality condition \eqref{E24} implies that $\bar u(x) = \alpha$, and hence $S(\bar u)(x) = \frac{1}{\tichonov}\bar y(x)\bar\varphi(x) <\alpha-\frac{\tau}{\tichonov}$. Using Lemma \ref{L4.2}, we have that
  \[S(u)(x) < S(\bar u)(x) + L_S \delta <\alpha-\frac{\tau}{\tichonov} +L_S \frac{\tau}{\tichonov L_S} = \alpha,\]
  and $x\in \mathbb{A}_u$ by definition of $\mathbb{A}_u$. The case $\tichonov\bar u(x)- \bar y(x)\bar\varphi(x) < \tau$ is treated in the same way using that, in this case $\bar u(x) = \beta$.

Therefore $(Q)$ has a unique solution $v\in L^2(\mathbb{I}_u)$.
Since $z_{u,\ci v}, \eta_{u,\ci v}\in L^s(\Gamma)$, \eqref{E28} implies that $v\in L^s(\mathbb{I}_u)$ and, consequently, $v$ is the unique solution of the equation $M_u v=w$ in $L^s(\Gamma)$.

To prove the uniform boundedness of $M_u^{-1}$ we proceed in two steps.

{\em Step 1.} Let us prove that there exists $C_2>0$ such that
\[\|v\|_{L^2(\Gamma)}\leq C_2 \|w\|_{L^s(\Gamma)}.\]
Since $\ci v\in T^\tau_{\bar u}$, we use the second order condition \eqref{E27}, the expression for the second derivative of $J$ obtained in Corollary \ref{C3.4}, equation \eqref{E28}  to obtain
\begin{align*}
  \frac{\kappa}{2}\| \ci v\|_{L^2(\Gamma)}^2 \leq  &
  J''(u) (\ci v)^2 =
  \int_\Gamma (\nu \ci v -\left(\varphi_u z_{u,\ci v}+ y_u \eta_{u,\ci v}\right) ) \ci v\dx \\
  =  & \int_\Gamma (\nu w + \left(\varphi_u z_{u,\ca w}+ y_u \eta_{u,\ca w}\right) ) \ci v\dx
\end{align*}
On the active se we have that $\ca w= \ca v$, so we can write
\[\nu\|\ca v\|^2_{L^2(\Gamma)} = \nu\int_\Gamma w \ca v\dx.\]
Therefore, adding the previous inequalities and applying Lemmas \ref{LA.2}, \ref{LA.5}, \ref{LA.7}, and \ref{LA.9} we obtain.
\begin{align*}
  \min\{\frac{\kappa}{2},\nu\}\| v\|^2_{L^2(\Gamma)} \leq & \nu \int_\Gamma  w v\dx+\int_\Gamma \left(\varphi_u z_{u,\ca w}+ y_u \eta_{u,\ca w}\right) ) \ci v\dx \\
  \leq  & |\Gamma|^{\frac{s-2}{2s}}[\nu + \sqrt{|\Gamma|}C_\Gamma (K_\infty^* C_G + K_\infty C_\Phi)] \|w\|_{L^s(\Gamma)} \|v\|_{L^2(\Gamma)},
\end{align*}
and we can take $C_2 =|\Gamma|^{\frac{s-2}{2s}}[\nu + \sqrt{|\Gamma|}C_\Gamma(K_\infty^* C_G + K_\infty C_\Phi)]/\min\{\frac{\kappa}{2},\nu\}$.

{\em Step 2.}  Finally, we prove that if $\dimension = 3$, then there exists $C>0$ such that
\[
\|v\|_{L^s(\Gamma)}\leq C \|w\|_{L^s(\Gamma)}.
\]
First, we use Lemmas \ref{LA.2}, \ref{LA.5}, \ref{LA.7} and \ref{LA.9}, and the boundedness of $M_u^{-1}$ in $L^2(\Gamma)$ and the estimate established in {\em Step 1}: for $C_3 =\sqrt{|\Gamma|}C_\Gamma(K_\infty^*C_G + K_\infty C_\Phi)$ we obtain
\[
\|\varphi_u z_{u,\ci v}+ y_u \eta_{u,\ci v}\|_{L^s(\Gamma)} \leq C_3 \|v\|_{L^2(\Gamma)} \leq C_3C_2  \|w\|_{L^s(\Gamma)}.
\]
Then, using \eqref{E28} and, once again Lemmas \ref{LA.2}, \ref{LA.5}, \ref{LA.7}, and \ref{LA.9} we have that
\[\|\ci v\|_{L^s(\Gamma)} \leq \left(1+\frac{C_3(1+C_2)}{\nu}\right) \|w\|_{L^s(\Gamma)}.\]
Since $\ca w= \ca v$, we conclude that $\|v\|_{L^s(\Gamma)}\leq C \|w\|_{L^s(\Gamma)}$, where $C =\max\{1, 1+\frac{C_3(1+C_2)}{\nu}\}$.
\end{proof}

Algorithm  \ref{Alg2} implements the semismooth Newton method to solve \Pb.
As a straightforward consequence of \cite[Theorem 3.13]{Ulbrich11}, Corollary \ref{C4.3}, and Theorem \ref{T4.4} we conclude the convergence of this algorithm.

\LinesNumbered
\begin{algorithm2e}[h!]
\caption{Semismooth Newton method for \Pb.}\label{Alg2}
\DontPrintSemicolon
Initialize. Choose $u_0\in \A$. Set $j=0$.\;
\For{$j \geq 0$}{
Compute $y_j = G(u_j)$\label{line3alg2}
\;
Compute $\varphi_j = \Phi(u_j)$
\;
Set $\mathbb{A}_j =\mathbb{A}^{\umax}_j\cup\mathbb{A}^{\umin}_j$ and  $\mathbb{I}_j=\Gamma\setminus\mathbb{A}_j$, where
\[\mathbb{A}^{\umax}_j = \{x\in\Gamma: y_j(x)\varphi_j(x)\geq \tichonov\umax\},\]
\[\mathbb{A}^{\umin}_j = \{x\in\Gamma: y_j(x)\varphi_j(x)\leq \tichonov\umin\}\]
\;
Set $w_j(x) = -F(u_j)(x)$:
\[w_j(x) =  \left\{\begin{array}{ll}
                              -u_j(x) +\umax & \text{ if }x \in \mathbb{A}^{\umax}_j \\
                              -u_j(x) + \frac{1}{\tichonov}\varphi_j(x)y_j(x) & \text{ if }x \in \mathbb{I}_j \\
                              -u_j(x) +\umin & \text{ if }x \in \mathbb{A}^{\umin}_j
                            \end{array}\right.
\]
\;
Compute $z_j = z_{u_j,\chi_{{}_{\mathbb{A}_j}}w_{j}}$ and $\eta_j = \eta_{u_j,\chi_{{}_{\mathbb{A}_j}}w_{j}}$\label{line8} \;
Solve the quadratic problem\label{line9}
\[
\mathrm{(}Q_j\mathrm{)}\quad \min_{v\in L^2(\mathbb{I}_j)} \mathbb{J}_j(v) := \frac{1}{2\tichonov} J''(u_j)(\chi_{{}_{\mathbb{I}_j}} v)^2
-\int_{\mathbb{I}_j}( w_j +\frac{1}{\tichonov}[z_j\varphi_j+y_j\eta_j]) v\dsigmax
\]
Name $v_{\mathbb{I}_j}$ its solution.\;
Set $u_{j+1}=u_j+\chi_{{}_{\mathbb{A}_j}}w_{j} + \chi_{{}_{\mathbb{I}_j}} v_{\mathbb{I}_j}$ and $j=j+1$.\;
}
\end{algorithm2e}

\begin{corollary}\label{C4.5}
Let $(\bar u,\bar y,\bar\varphi) \in \uad \times Y^2$ satisfy the first order optimality conditions \eqref{E22}--\eqref{E24}, the strict complementarity condition $|\Sigma_{\bar u}| = 0$, and the second order sufficient optimality condition $J''(\bar u)v^2 > 0$ for every $v \in C_{\bar u}\setminus\{0\}$. Then, there exists $\delta > 0$ such that for all $u_0 \in B_\delta(\bar u)$ the sequence $\{u_j\}$ generated by Algorithm \ref{Alg2} is contained in the ball $B_\delta(\bar u)$ and converges superlinearly to $\bar u$.
\end{corollary}
The radius of the basin of attraction $\delta$ depends on parameters related to the continuity properties of the involved functionals and its derivatives, the second order condition and the neighborhood in $L^s(\Gamma)$ for which the state equation is meaningful.

\section{A numerical example and some computational considerations}
Consider $\Omega=(0,1)^3$, $Ay=-\Delta y+y$, $a(x,y)= y^3-\sin(2\pi x_1)\sin(\pi x_2)\cos(3\pi x_3)$, $g\equiv 0$, $L(x,y) = 0.5(y-y_d(x))^{2}$, with $y_d(x) = -512\prod_{i=1}^3 x_i(1-x_i)$, $\tichonov = 0.01$, $\umin=0$, and $\umax = 1$.
We solve a finite element discretization of \Pb. Continuous piecewise linear functions are used for the state, the adjoint state, and the control. The Tichonov regularization term is discretized using the lumped mass matrix. In this way, the optimization algorithm for the discrete problem is exactly the discrete version of Algorithm \ref{Alg2}.

The convergence history for $u_0=0$ is summarized in tables \ref{T1} and \ref{T2} for different mesh sizes. The expected superlinear convergence can
be seen in the relative errors between consecutive iterations, denoted $\delta_j$.
We also remark the mesh-independence of the convergence history, which is to be expected since we have obtained our results in the infinite-dimensional setting.

At each iteration we have to solve a nonlinear equation to compute $y_j$ and solve an unconstrained quadratic problem to compute $v_{\mathbb{I}_j}$. We use Newton's method for the first task and the conjugate gradient method for the second one. Notice that $\mathbb{J}_j(v) = \frac{1}{2}(v,A_j v)_{L^2(\mathbb{I}_j)} - (b_j,v)_{L^2(\mathbb{I}_j)}$, where $b_j = \chi_{{}_{\mathbb{I}_j}} ( w_{j} +\frac{1}{\tichonov}[z_j\varphi_j+y_j\eta_{j}])$ and, for any $v\in L^2(\mathbb{I}_j)$, \[A_j v = \chi_{{}_{\mathbb{I}_j}} \left(v+\frac{1}{\tichonov}[z_{u_j,\chi_{{}_{\mathbb{I}_j}} v}\varphi_j+\eta_{u_j,\chi_{{}_{\mathbb{I}_j}} v} y_j]\right);\]
see eqs. \eqref{E13} and \eqref{E21}

We include in the tables the number of Newton iterations used to solve the nonlinear equation at each iteration. Each of these requires the factorization of the finite element matrix, and this number is a good measure of the
global
complexity of the method. In contrast, each of the conjugate gradient iterations used to solve $\mathrm{(}Q_j\mathrm{)}$
requires the solution of two linear systems, but the matrix has been previously factorized in the last step of the nonlinear solve.

\begin{table}[h!]
  \centering
  \begin{tabular}{ccccc}
    $j$ &       $J(u_j)$            &   $\delta_j$   &  $\sharp$Newton & $\sharp$CG  \\ \hline
  0&  4.7607853276096295e+00 &   7.3e-01  &   3 &  17  \\
  1&  4.7590621154705985e+00 &   5.3e-01  &   3 &  12  \\
  2&  4.7588905662088630e+00 &   1.1e-01  &   3 &  12  \\
  3&  4.7588301468521248e+00 &   3.7e-04  &   3 &  12  \\
  4&  4.7588301456859448e+00 &   7.9e-08  &   2 &  12  \\
  5&  4.7588301456859456e+00 &   3.7e-15  &   2 &  12  \\
  6 & 4.7588301456859456e+00 &            &   1 &     \end{tabular}
  \caption{Solution of  \Pb for $h=2^{-4}$. }\label{T1}
\end{table}

\begin{table}[h!]
  \centering
  \begin{tabular}{ccccc}
    $j$ &       $J(u_j)$            &   $\delta_j$   &  $\sharp$Newton & $\sharp$CG  \\ \hline
  0&  4.8308890801571112e+00 &   7.9e-01  &   3 &  16  \\
  1&  4.8290362150750905e+00 &   5.8e-01  &   3 &  11  \\
  2&  4.8288131518545896e+00 &   1.3e-01  &   3 &  12  \\
  3&  4.8287240470263058e+00 &   7.3e-04  &   3 &  11  \\
  4&  4.8287240439741863e+00 &   4.7e-06  &   2 &  11  \\
  5&  4.8287240439742973e+00 &   6.3e-14  &   2 &  11  \\
  6 & 4.8287240439742973e+00 &            &   1 &
  \end{tabular}
  \caption{Solution of  \Pb for  $h=2^{-5}$. }\label{T2}
\end{table}


\providecommand{\bysame}{\leavevmode\hbox to3em{\hrulefill}\thinspace}
\providecommand{\MR}{\relax\ifhmode\unskip\space\fi MR }
\providecommand{\MRhref}[2]{%
  \href{http://www.ams.org/mathscinet-getitem?mr=#1}{#2}
}
\providecommand{\href}[2]{#2}

\appendix

\section{Proofs of some auxiliary results}\label{Apx}

\begin{lemma}\label{LA.1}
  For every $u\in\A$ and every $y\in L^\infty(\Omega)$,
 \[
     \mathfrak{a}(z,z) + \int_\Omega \displaystyle\frac{\partial a}{\partial y}(x,y)z^2\dx + \int_\Gamma  u z^2\dx
     \geq
      \frac{\lambda_A}{2} \|z\|^2_{H^1(\Omega)} \ \forall z\in H^1(\Omega).
  \]
\end{lemma}
\begin{proof}
  From the construction of $\A$, we know that there exists $\bar u\in \A_0$ such that $\|u-\bar u\|_{L^s(\Gamma)}<\varepsilon_{\bar u}$, with $\varepsilon_{\bar u} < \frac{1}{2}\lambda_A/(|\Gamma|^{\frac{s-2}{s}}C_\Gamma^2)$.

  Using assumptions \ref{A2.1} and \ref{A2.2}, we have that
   \begin{align*}
         \mathfrak{a}(z,z) & + \int_\Omega \displaystyle\frac{\partial a}{\partial y}(x,y)z^2\dx + \int_\Gamma  u z^2\dx \\
         & =
         \mathfrak{a}(z,z) + \int_\Omega \displaystyle\frac{\partial a}{\partial y}(x,y)z^2\dx + \int_\Gamma \bar u z^2\dx
     + \int_\Gamma ( u-\bar u)z^2\dx \\
      & \geq  \lambda_A \|z\|^2_{H^1(\Omega)} -  \int_\Gamma | u-\bar u|z^2\dx.
  \end{align*}
  Using the Cauchy-Schwarz inequality, \eqref{E05}, and the properties of $\varepsilon_{\bar u}$, we have
  \[\int_\Gamma | u-\bar u|z^2\dx\leq
  \|\bar u-u\|_{L^s(\Gamma)} \|z\|^2_{L^4(\Gamma)} |\Gamma|^{\frac{s-2}{s}}
  \leq
  \varepsilon_{\bar u}|\Gamma|^{\frac{s-2}{s}}
     C_\Gamma^2\|z\|^2_{H^1(\Omega)}\leq \frac{\lambda_A}{2}\|z\|^2_{H^1(\Omega)},\]
     and the proof is complete.
\end{proof}
\begin{lemma}\label{LA.2}
  There exist constants $C'$, $M'_\infty$, $K_\infty$ and $C'_{\mu,\infty}$ such that, for every $u\in\A$
  \begin{align}
  \|y_u\|_{H^1(\Omega)}\leq & C'(\|a(\cdot,0)\|_{L^p(\Omega)}+\|g\|_{L^q(\Omega)}),\label{EA.1}\\
  \|y_u\|_{L^\infty(\Omega)}\leq & M'_\infty(\|a(\cdot,0)\|_{L^p(\Omega)}+\|g\|_{L^q(\Omega)}) = :K_\infty,\label{EA.2}\\
  \|y_u\|_{C^{0,\mu}(\bar\Omega)} \leq & C'_{\mu,\infty}(\|a(\cdot,0)\|_{L^p(\Omega)}+\|u\|_{L^s(\Gamma)}+\|g\|_{L^q(\Gamma)}).
  \label{EA.3}
  \end{align}
\end{lemma}
\begin{proof}
  Given $u\in\A$, we take  $\bar u\in \A_0$ such that $\|u-\bar u\|_{L^s(\Gamma)}<\varepsilon_{\bar u}$.
   Denote $z =y_{\bar u}-y_u\in Y$. Subtracting the equations satisfied by $y_{\bar u}$ and $y_u$ and using the mean value theorem, we obtain
   \begin{equation}\label{EA.M1}\left\{
   \begin{array}{l}
     A z + \displaystyle\frac{\partial a}{\partial y}(x,y_\theta)z =0\mbox{ in }\Omega, \\
     \partial_{\conormal_A} z + uz = (u-\bar u) y_{\bar u}\text{ on }\Gamma,
   \end{array}
   \right.
   \end{equation}
   where $y_\theta = y_u+\theta(y_{\bar u}-y_u)$ for a measurable function  $0\leq \theta\leq 1$.
   With the help of \eqref{E08}, we notice that
   \[\| (u-\bar u) y_{\bar u}\|_{L^s(\Gamma)} \leq \varepsilon_{\bar u} M_\infty(\|a(\cdot,0)\|_{L^p(\Omega)}+
   \|g\|_{L^q(\Gamma)}).\]
Hence, thanks to Lemma \ref{LA.1}, applying the methods of \cite[Theorem 4.1]{Stampacchia65}, we infer that
\[\|z\|_{H^1(\Omega)}+\|z\|_{L^\infty(\Omega)}\leq C_1(\|a(\cdot,0)\|_{L^p(\Omega)}+
   \|g\|_{L^q(\Gamma)}).\]
Then, using that
\[\|y_u\|_{H^1(\Omega)}+\|y_u\|_{L^\infty(\Omega)} \leq \|z\|_{H^1(\Omega)}+\|z\|_{L^\infty(\Omega)} +
\|y_{\bar u}\|_{H^1(\Omega)}+\|y_{\bar u}\|_{L^\infty(\Omega)},\]
the estimates \eqref{EA.1} and \eqref{EA.2} follow from this inequality, the above estimate for $z$, and Theorem \ref{T2.3}.
     Finally, \eqref{EA.3} is obtained using the same technique as for \eqref{E09}, but using \eqref{EA.2} instead of \eqref{E08}.
\end{proof}
\begin{lemma}\label{LA.3}
  The solution mapping $G:\A\to Y$ is locally Lipschitz: for every $u\in\A$, there exist $\delta_u>0$ and $L_{G,u}>0$ such that
  \[\|y_{u_1}-y_{u_2}\|_Y \leq L_{G,u}\|u_1-u_2\|_{L^s(\Gamma)}\ \forall u_1,u_2\in B_{\delta_u}(u).\]
\end{lemma}
\begin{proof}
Since $G:\A\to Y$ is of class $C^1$, the mapping $DG:\A\longrightarrow\mathcal{L}(L^s(\Gamma),Y)$ is continuous. Therefore, given $u\in\A$ there exist $\delta_u>0$ and $L_{G,u}$ such that $B_{\delta_u}(u)\subset\A$ and $\|DG(\hat u)\|_{\mathcal{L}(L^s(\Gamma),Y)}\leq L_{G,u}$ for evey $\hat u\in B_{\delta_{u}}(u)$. The Lipschit property on this ball is a straightforward consequence of the generalized mean value theorem.
\end{proof}

\begin{lemma}\label{LA.4}
  For every $u\in\A$, every $f\in L^2(\Omega)$ and every $v\in L^2(\Gamma)$, the equation
  \[\left\{\begin{array}{l} A\zeta + \displaystyle \frac{\partial a}{\partial y}(x,y_u)\zeta = f\ \  \mbox{in } \Omega,\vspace{2mm}\\  \partial_{\conormal_A} \zeta + u\zeta  = v \ \ \mbox{on }\Gamma, \end{array}\right.
  \]
 has a unique solution $\zeta\in H^1(\Omega)$ and there exists a constant $C_{A}>0$ independent of $u$ and $v$ such that
  \[\|\zeta\|_{H^1(\Omega)}\leq C_{A} (\|f\|_{L^2(\Omega)} + \|v\|_{L^2(\Gamma)}).\]
  If we replace the operator $A$ by $A^*$ in the previous equation, the statement stays true and the inequality holds with the same constant $C_A$.
\end{lemma}
\begin{proof}
  Take $v\in L^2(\Gamma)$. There exists a sequence $\{v_k\}_k\subset L^s(\Gamma)$ such that $v_k\to v$ in $L^2(\Gamma)$. By Theorem \ref{T2.4}, there exists $\zeta_k\in Y$ such that
  \begin{equation}
  \label{EA.4}\mathfrak{a}(\zeta_k,\phi)+\int_\Omega \frac{\partial a}{\partial y}(x,y_u)\zeta_k\phi\dx + \int_\Gamma u \zeta_k\phi \dx=\int_\Omega f \phi\dx +
  \int_\Gamma v_k \phi\dx\ \forall\phi\in H^1(\Omega)
  \end{equation}
Testing the variational formulation for $\phi=\zeta_k$ and using Lemma \ref{LA.1} and assumptions \ref{A2.1} and \ref{A2.2} we have
  \begin{align*}
    \frac{\lambda_A}{2}\|\zeta_k\|^2_{H^1(\Omega)} \leq  &
    \mathfrak{a}(\zeta_k,\zeta_k)+\int_\Omega \frac{\partial a}{\partial y}(x,y_u)\zeta_k^2\dx + \int_\Gamma  u \zeta_k^2\dx
    =
      \int_\Omega f \phi\dx +\int_\Gamma v_k\zeta_k\dx \\
    \leq & \|f\|_{L^2(\Omega)}\|\zeta_k\|_{L^2(\Omega)}+ \|v_k\|_{L^2(\Gamma)} \|\zeta_k\|_{L^2(\Gamma)}
    \\
    \leq
    &
    (\|f\|_{L^2(\Omega)} + C_\Gamma |\Gamma|^{1/4} \|v_k\|_{L^2(\Gamma)}) \|\zeta_k\|_{H^1(\Omega)}
  \end{align*}
  Dividing by $\|\zeta_k\|_{H^1(\Omega)}$, we get
  \begin{equation}\label{EA.M2}
  \|\zeta_k\|_{H^1(\Omega)}\leq C_{A} (\|f\|_{L^2(\Omega)} +\|v_k\|_{L^2(\Gamma)}),\end{equation}
  where
  \[C_{A} = \frac{2}{\lambda_A} \max\{1,|\Gamma|^{1/4} C_\Gamma\}.\]
  Since the sequence $\{v_k\}_k$ is bounded in $L^2(\Gamma)$, then $\{\zeta_k\}_k$ is bounded in $H^1(\Omega)$. Thus, we can extract a subsequence, denoted in the same way, such that $\zeta_k\rightharpoonup \zeta$ weakly in $H^1(\Omega)$. Taking limits in \eqref{EA.4}, and \eqref{EA.M2}, we get that $\zeta$ solves the variational formulation of the equation and the claimed estimate is satisfied.
\end{proof}
\begin{lemma}\label{LA.5}
  For every $u\in\A$ and every $v\in L^2(\Gamma)$, the equation \eqref{E13} has a unique solution $z_{u,v}\in H^1(\Omega)$ and there exists a constant $C_G>0$ independent of $u$ and $v$ such that
  \[\|z_{u,v}\|_{H^1(\Omega)}\leq C_G \|v\|_{L^2(\Gamma)}.\]
\end{lemma}
\begin{proof}
The result follows from Lemma \ref{LA.4} taking into account Lemma \ref{LA.2} and using $C_G = K_\infty C_{A}$.
  \end{proof}
\begin{lemma}\label{LA.6}
 For every $u\in \A$ there exists $L_{G',u}>0$ such that
  \[\|z_{u_1,v}-z_{u_2,v}\|_{H^1(\Omega)}\leq L_{G',u} \|u_1-u_2\|_{L^s(\Gamma)} \|v\|_{L^2(\Gamma)}\ \forall v\in L^2(\Gamma)\ \forall u_1,u_2\in B_{\delta_u}(u),\]
where $\delta_u>0$ is the one introduced in Lemma \ref{LA.3}.
 \end{lemma}
\begin{proof}
  Denote $\zeta = z_{u_1,v}-z_{u_2,v}\in H^1(\Omega)$. This function satisfies
  \[\left\{
   \begin{array}{l}
     A \zeta + \displaystyle\frac{\partial a}{\partial y}(x,y_{u_1})\zeta = \left(\displaystyle\frac{\partial a}{\partial y}(x,y_{u_2}) -\displaystyle\frac{\partial a}{\partial y}(x,y_{u_1}) \right) z_{u_2,v} \mbox{ in }\Omega \\
     \partial_{\conormal_A} \zeta + u_1 \zeta = (u_2- u_1) z_{u_2,v} + v( y_{u_2}-y_{u_1}) \text{ on }\Gamma.
   \end{array}
   \right.
   \]
   We estimate the right hand side of the above equation:
   \begin{align*}
     \|\left(\frac{\partial a}{\partial y}(x,y_{u_2}) \right. & \left.
      -\frac{\partial a}{\partial y}(x,y_{u_1}) \right) z_{u_2,v}\|_{L^2(\Omega)}  \\
      =  & \|\frac{\partial^2 a}{\partial y^2}(x,y_\theta)(y_{u_2}-y_{u_1}) z_{u_2,v}\|_{L^2(\Omega)}  \\
      \leq  & \|\frac{\partial^2 a}{\partial y^2}(x,y_\theta)\|_{L^\infty(\Omega)}\|(y_{u_2}-y_{u_1})\|_{L^\infty(\Omega)} \|z_{u_2,v}\|_{L^2(\Omega)} \\
      \leq & C_{a,K_\infty} L_{G,u} C_G \|u_1-u_2\|_{L^s(\Gamma)} \|v\|_{L^2(\Gamma)}.
   \end{align*}
   Now, we estimate the boundary terms. For the first term we get with Lemma \ref{LA.5}
   \begin{align*}
     \|(u_2- u_1) z_{u_2,v} \|_{L^2(\Gamma)} \leq
     &
     \|u_1-u_2\|_{L^s(\Gamma)}  \|z_{u_2,v}\|_{L^4(\Gamma)} |\Gamma|^{\frac{s-2}{s}}\\
     \leq &
     C_\Gamma C_G |\Gamma|^{\frac{s-2}{s}}\|u_1-u_2\|_{L^s(\Gamma)}\|v\|_{L^2(\Gamma)}.
   \end{align*}
   For the second term we have
   \begin{align*}
     \|v( y_{u_2}-y_{u_1})\|_{L^2(\Gamma)} \leq
     \| y_{u_2}-y_{u_1}\|_{L^\infty(\Omega)} \|v\|_{L^2(\Gamma)}
     \leq
     L_{G,u} \|u_1-u_2\|_{L^s(\Gamma)} \|v\|_{L^2(\Gamma)}.
   \end{align*}
The proof concludes by straightforward application of Lemma \ref{LA.4} taking
\[L_{G',u} = C_{A}\left(L_{G,u} ( C_{a,K_\infty}  C_G +1 ) + C_\Gamma C_G |\Gamma|^{\frac{s-2}{s}} \right).\]
\end{proof}
\begin{lemma}\label{LA.7}
  For every $u\in\A$, $\|\varphi_u\|_{L^\infty(\Omega)}\leq K^*_\infty$, where $K^*_\infty$ is independent of $u$.
\end{lemma}
\begin{proof}
  Applying Lemma \ref{LA.2} to the adjoint state equation and using that $\|y_u\|_{L^\infty(\Omega)}\leq K_\infty$ and Assumption \ref{A3.1}, we obtain the existence of a constant $M^*_\infty>0$ such that
  \[\|\varphi_u\|_{L^\infty(\Omega)} \leq M^*_\infty C_{L,K_\infty}=:K_\infty^*.\]
\end{proof}
\begin{lemma}\label{LA.8}
  The mapping $\Phi:\A\to Y$ is locally Lipschitz. for every $u\in\A$, there exist $\delta_u>0$ and $L_{\Phi,u}>0$ such that
  \[\|\varphi_{u_1}-\varphi_{u_2}\|_Y \leq L_{\Phi,u}\|u_1-u_2\|_{L^s(\Gamma)}\ \forall u_1,u_2\in B_{\delta_u}(u).\]
\end{lemma}
\begin{proof}Since $\Phi:\A\longrightarrow Y$ is $C^1$, arguing as in the proof of Lemma \ref{LA.3}, the statement follows.
\end{proof}
\begin{lemma}\label{LA.9}
    For every $u\in\A$ and every $v\in L^2(\Gamma)$, the equation \eqref{E21} has a unique solution $\eta_{u,v}\in H^1(\Omega)$ and there exists a constant $C_\Phi>0$ independent of $u$ and $v$ such that
  \[\|\eta_{u,v}\|_{H^1(\Omega)}\leq C_\Phi \|v\|_{L^2(\Gamma)}.\]
\end{lemma}
\begin{proof}
Working as in the proof of Lemma \ref{LA.6}, we have that
\begin{align*}
  &  \|
    \Big[\frac{\partial^2 L}{\partial y^2}(x,y_u) - \varphi_u  \frac{\partial^2 a}{\partial y^2}(x,y_u)\Big] z_{u,v}
    \|_{L^2(\Omega)} \leq
    \|
    \frac{\partial^2 L}{\partial y^2}(x,y_u) - \varphi_u  \frac{\partial^2 a}{\partial y^2}(x,y_u)\|_{L^\infty(\Omega)}
    \| z_{u,v}    \|_{L^2(\Omega)}
     \\
   & \leq  (C_{L,K_\infty} + K^*_\infty C_{a,K_\infty})
    \|z_{u,v}\|_{H^1(\Omega)} \leq (C_{L,K_\infty} + K^*_\infty C_{a,K_\infty})
    C_G \|v\|_{L^2(\Gamma)},
\end{align*}
and that $\|v\varphi_u\|_{L^2(\Gamma)}\leq K^*_\infty \|v\|_{L^2(\Gamma)}$. Therefore, the result is an immediate consequence of Lemma \ref{LA.2} taking
\[C_\Phi = C_{A} \left((C_{L,K_\infty} + K^*_\infty C_{a,K_\infty})    C_G + K^*_\infty\right).\]
\end{proof}
\begin{lemma}\label{LA.10}
    For every $u\in \A$ and for every $\varepsilon >0 $ there exists $\rho^*_{\varepsilon, u}>0$ such that
  \[\|\eta_{u_1,v}-\eta_{ u_2,v}\|_{H^1(\Omega)}\leq \varepsilon \|v\|_{L^2(\Gamma)}\ \forall v\in L^2(\Gamma)\ \forall u_1,u_2\in B_{\rho^*_{\varepsilon, u}}(u).\]
\end{lemma}
\begin{proof}Take  $\rho_0 = \delta_{u}$, where $\delta_u$ is the minimum of the ones introduced in Lemmas \ref{LA.3} and \ref{LA.8} and assume that $u_1,u_2\in B_{\rho_0}(u)$.
   Define $\zeta = \eta_{u_1,v}-\eta_{ u_2,v}$. This function satisfies
    \[\left\{
   \begin{array}{l}
     A^* \zeta + \displaystyle\frac{\partial a}{\partial y}(x,y_{u_1})\zeta = \left(\displaystyle\frac{\partial a}{\partial y}(x,y_{u_2}) -\displaystyle\frac{\partial a}{\partial y}(x,y_{u_1}) \right) z_{u_2,v} \\
 \displaystyle   + \Big[\frac{\partial^2 L}{\partial y^2}(x,y_{u_1}) - \varphi_{u_1}  \frac{\partial^2 a}{\partial y^2}(x,y_{u_1})\Big] z_{u_1,v}
    - \Big[\frac{\partial^2 L}{\partial y^2}(x,y_{u_2}) - \varphi_{u_2}  \frac{\partial^2 a}{\partial y^2}(x,y_{u_2})\Big] z_{u_2,v} \mbox{ in }\Omega, \\
     \partial_{\conormal_{A^*}} \zeta + u_1 \zeta = (u_2- u_1) \eta_{u_2,v} \text{ on }\Gamma.
   \end{array}
   \right.
   \]

   We are going to apply Lemma \ref{LA.4}. To this end it is enough to estimate the right hand side of the equation in $L^2(\Omega)\times L^2(\Gamma)$ by $\varepsilon \|v\|_{L^2(\Gamma)}$.

   \begin{align*}
     &  \left[\frac{\partial^2 L}{\partial y^2}(x,y_{u_1}) - \varphi_{u_1}  \frac{\partial^2 a}{\partial y^2}(x,y_{u_1})\right] z_{u_1,v}
    - \left[\frac{\partial^2 L}{\partial y^2}(x,y_{u_2}) - \varphi_{u_2}  \frac{\partial^2 a}{\partial y^2}(x,y_{u_2})\right] z_{u_2,v}  \\
  = &     \left[
     \frac{\partial^2 L}{\partial y^2}(x,y_{u_1})
      -\frac{\partial^2 L}{\partial y^2}(x,y_{u_2})
      \right] z_{u_1,v}
      \\
    &  +\varphi_{u_1}\left[\frac{\partial^2 a}{\partial y^2}(x,y_{u_2})-\frac{\partial^2 a}{\partial y^2}(x,y_{u_1})\right]z_{u_1,v}
       +(\varphi_{u_2}-\varphi_{u_1})  \frac{\partial^2 a}{\partial y^2}(x,y_{u_2})
      z_{u_1,v}\\
    & +\left[\frac{\partial^2 L}{\partial y^2}(x,y_{u_2}) - \varphi_{u_2}  \frac{\partial^2 a}{\partial y^2}(x,y_{u_2})\right] (z_{u_1,v}-z_{u_2,v}).
   \end{align*}

{\em Estimation of the first term}. Consider $\varepsilon_1 = \frac{\varepsilon}{6 C_G C_{A}}$. From Assumption \ref{A3.1} we know that there exists $\tilde \rho_1>0$ such that if
\begin{equation}\label{EA.5}
\|y_{u_1}-y_{u_2}\|_{L^\infty(\Omega)} < \tilde\rho_1,
\end{equation}
 then
\begin{equation}\label{EA.6}
\|
     \frac{\partial^2 L}{\partial y^2}(x,y_{u_1})
      -\frac{\partial^2 L}{\partial y^2}(x,y_{u_2})
      \|_{L^\infty(\Omega)} < \varepsilon_1.
\end{equation}
From Lemma \ref{LA.3} we infer the existence of $\delta_{u,1}>0$ such that
\[
\|y_{u_1}-y_{u_2}\|_{L^\infty(\Omega)}\leq L_{G,u} \|u_1-u_2\|_{L^s(\Gamma)}
\quad \forall u_1,u_2\in B_{\delta_{u,1}}(u).
\]
Define $\rho_1 = \min\{\frac{\tilde\rho_1}{2L_{G,u}},\delta_{u,1}\}.$
Hence, if $u_1,u_2\in B_{\rho_1}(u)$ we have that \eqref{EA.5} holds and consequently, so does \eqref{EA.6}.
Using this,
we deduce, with the help of Assumption \ref{A3.1} and Lemma \ref{LA.5} that
\begin{align*}
  \|   \left[
     \frac{\partial^2 L}{\partial y^2}(x,y_{u_1}) \right.&
    \left.  -\frac{\partial^2 L}{\partial y^2}(x,y_{u_2})
      \right] z_{u_1,v} \|_{L^2(\Omega)} \\
 \leq &     \|
     \frac{\partial^2 L}{\partial y^2}(x,y_{u_1})
      -\frac{\partial^2 L}{\partial y^2}(x,y_{u_2})
      \|_{L^\infty(\Omega)}
      \| z_{u_1,v} \|_{L^2(\Omega)}
      \leq
       \frac{\varepsilon}{6C_{A}} \|v\|_{L^2(\Gamma)}.
\end{align*}

{\em Estimation of the second term}. Consider $\varepsilon_2 = \frac{\varepsilon}{6 K^*_\infty C_G C_{A}}$. From Assumption \ref{A2.2} we know that there exists $\tilde\rho_2>0$ such that, if
\begin{equation}\label{EA.7}
\|y_{u_1}-y_{u_2}\|_{L^\infty(\Omega)} < \tilde\rho_2,
\end{equation}
 then
\begin{equation}\label{EA.8}
\|
     \frac{\partial^2 a}{\partial y^2}(x,y_{u_1})
      -\frac{\partial^2 a}{\partial y^2}(x,y_{u_2})
      \|_{L^\infty(\Omega)} < \varepsilon_2.
\end{equation}
With Lemma \ref{LA.3} we deduce the existence of $\delta_{u,2}>0$ such that
\[
\|y_{u_1}-y_{u_2}\|_{L^\infty(\Omega)}\leq L_{G,u} \|u_1-u_2\|_{L^s(\Gamma)}
\quad \forall u_1,u_2\in B_{\delta_{u,2}}(u).
\]
Define $\rho_2 = \min\{\frac{\tilde\rho_2}{2L_{G,u}},\delta_{u,2}\}$. Hence, if $u_1,u_2\in B_{\rho_2}(u)$ we have that \eqref{EA.7} holds and so does \eqref{EA.8}. Using this,
we deduce with the help of Assumption \ref{A2.2}, Lemma \ref{LA.7}, and Lemma \ref{LA.5} that
\begin{align*}
  \|   \varphi_{u_1} & \left[
     \frac{\partial^2 a}{\partial y^2}(x,y_{u_1})
      -\frac{\partial^2 a}{\partial y^2}(x,y_{u_2})
      \right] z_{u_1,v} \|_{L^2(\Omega)} \\
 \leq & \|    \varphi_{u_1} \|_{L^\infty(\Omega)}
 \|   \frac{\partial^2 a}{\partial y^2}(x,y_{u_1})
      -\frac{\partial^2 a}{\partial y^2}(x,y_{u_2})
      \|_{L^\infty(\Omega)}
      \| z_{u_1,v} \|_{L^2(\Omega)}
      \leq
       \frac{\varepsilon}{6C_{A}} \|v\|_{L^2(\Gamma)}.
\end{align*}

{\em Estimation of the third term}. If $u_1,u_2\in B_{\rho_3}(u)$ with $\rho_3 = \frac{\varepsilon}{12 L_{\Phi,u} C_{a,K_\infty} C_G C_{A}}$,
we can deduce, using Lemma \ref{LA.7}, assumption \ref{A2.2} and Lemma \ref{LA.5} that
\begin{align*}
  \|(\varphi_{u_2}-\varphi_{u_1})&  \frac{\partial^2 a}{\partial y^2}(x,y_{u_2})
      z_{u_1,v}\|_{L^2(\Omega)} \\
      \leq  &
      \|(\varphi_{u_2}-\varphi_{u_1}) \|_{L^\infty(\Omega)} \|\frac{\partial^2 a}{\partial y^2}(x,y_{u_2})\|_{L^\infty(\Omega)}
      \|z_{u_1,v}\|_{L^2(\Omega)}      \\
  \leq  & L_{\Phi,u}\|u_1-u_2\|_{L^s(\Gamma)} C_{a,K_\infty} C_G \|v\|_{L^2(\Gamma)}\leq
       \frac{\varepsilon}{6C_{A}} \|v\|_{L^2(\Gamma)}.
\end{align*}

{\em Estimation of the fourth term}. If $u_1,u_2\in B_{\rho_4}(u)$ with $\rho_4 = \frac{\varepsilon}{12 L_{G',u} (C_{L,K_\infty} + K^*_\infty C_{a,K_\infty})C_{A}}$, we deduce with the help of Lemmas \ref{LA.2} and \ref{LA.7}, assumptions \ref{A2.2} and \ref{A3.1}, and Lemma \ref{LA.6} that
\begin{align*}
  \| \left[\frac{\partial^2 L}{\partial y^2}(x,y_{u_2}) - \right. & \left. \varphi_{u_2}  \frac{\partial^2 a}{\partial y^2}(x,y_{u_2})\right] (z_{u_1,v}-z_{u_2,v}) \|_{L^2(\Omega)} \\
   \leq & \| \frac{\partial^2 L}{\partial y^2}(x,y_{u_2}) - \varphi_{u_2}  \frac{\partial^2 a}{\partial y^2}(x,y_{u_2})\|_{L^\infty(\Omega)} \|z_{u_1,v}-z_{u_2,v}\|_{L^2(\Omega)} \\
\leq &
   \left(C_{L,K_\infty} + K^*_\infty C_{a,K_\infty}\right) L_{G',u} \|u_1-u_2\|_{L^s(\Gamma)} \|v\|_{L^2(\Gamma)}\leq
       \frac{\varepsilon}{6C_{A}} \|v\|_{L^2(\Gamma)}.
\end{align*}

{\em Estimation of the fifth term}. If $u_1,u_2\in B_{\rho_5}(u)$ with $\rho_5 = \frac{\varepsilon}{12 C_{a,K_\infty} L_{G,u} C_G C_{A}}$, we have, using the mean value theorem, Assumption \ref{A2.2}, Lemmas \ref{LA.2}, \ref{LA.3} and \ref{LA.5}
\begin{align*}
 \| \left(\frac{\partial a}{\partial y}(x,y_{u_2}) \right.&\left. -\displaystyle\frac{\partial a}{\partial y}(x,y_{u_1}) \right) z_{u_2,v}\|_{L^2(\Omega)}  \leq
 \| \frac{\partial a}{\partial y}(x,y_{u_2}) -\displaystyle\frac{\partial a}{\partial y}(x,y_{u_1}) \|_{L^\infty(\Omega)} \| z_{u_2,v}\|_{L^2(\Omega)}
  \\
  \leq &
  \| \frac{\partial^2 a}{\partial y^2}(x,y_{\theta}) \|_{L^\infty(\Omega)} \|y_{u_1}-y_{u_2} \|_{L^\infty(\Omega)} C_G \| v\|_{L^2(\Gamma)}\\
  \leq & C_{a,K_\infty} L_{G,u} C_G \|u_1-u_2\|_{L^s(\Gamma)} \|v\|_{L^2(\Gamma)} \leq \frac{\varepsilon}{6C_{A}}\|v\|_{L^2(\Gamma)}.
\end{align*}
{\em Estimation of the boundary term}. If $u_1,u_2\in B_{\rho_6}(u)$ with $\rho_6 = \frac{\varepsilon}{12 C_\Gamma C_\Phi |\Gamma|^{\frac{s-2}{s}}C_{A}}$, then using \eqref{E05} and Lemma \ref{LA.9}, we obtain
\begin{align*}
\|(u_1-u_2)\eta_{u_2,v}\|_{L^2(\Gamma)} \leq & \|u_1-u_2\|_{L^s(\Gamma)} \|\eta_{u_2,v}\|_{L^4(\Gamma)}
|\Gamma|^{\frac{s-2}{s}}\\
\leq & \rho_6 C_\Gamma C_\Phi |\Gamma|^{\frac{s-2}{s}} \|v\|_{L^2(\Gamma)} =
\frac{\varepsilon}{6C_{A}}\|v\|_{L^2(\Gamma)}.
\end{align*}
The proof concludes taking  $\rho^*_{u,\varepsilon} = \min\{\rho_i,\,i=0,\ldots,6\}$ and applying Lemma \ref{LA.4}.
\end{proof}
\begin{lemma}\label{LA.11}
  For every $u\in\A$ and every $\varepsilon>0$ there exists $\rho_{u,\varepsilon} >0$ such that
\[|(J''(u_1)-J''(u_2))v^2|<\varepsilon \|v\|_{L^2(\Gamma)}^2\ \forall v\in L^2(\Gamma)\ \forall u_1,u_2\in B_{\rho_{u,\varepsilon}}(u).\]
\end{lemma}
\begin{proof}Define $\rho_0 = \rho^*_{u,\varepsilon}$, where  $\rho^*_{u,\varepsilon}$ is defined in Lemma \ref{LA.10} and take $u_1,u_2\in B_{\rho_0}(u)$.
  Using Corollary \ref{C3.4} we have that
  \begin{align*}
    |(J''(u_1)&-J''(u_2))v^2 | =
    \left|\int_\Gamma \left(\varphi_{u_1}z_{u_1,v}+y_{u_1}\eta_{u_1,v} -
    (\varphi_{u_2}z_{u_2,v}+y_{u_2}\eta_{u_2,v} ) \right) v\dx\right|
     \\
    \leq  &     \int_\Gamma \left|
    \varphi_{u_1}(z_{u_1,v}-z_{u_2,v})v\right| \dx+
    \int_\Gamma \left| (\varphi_{u_1}-\varphi_{u_2})z_{u_2,v}v\right| \dx\\
    &  +
    \int_\Gamma \left| y_{u_1}(\eta_{u_1,v} - \eta_{u_2,v}) v \right| \dx +
    \int_\Gamma \left| (y_{u_1}-y_{u_2})\eta_{u_2,v}  v\right| \dx \\
    = & I + II + III + IV
  \end{align*}
  Define $\rho_1 =  \frac{\varepsilon}{8|\Gamma|^{1/4}C_\Gamma K^*_\infty L_{G',u}}$. If $u_1,u_2\in B_{\rho_1}(u)$, using Cauchy inequality, \eqref{E05}, Lemma \ref{LA.6}, and Lemma \ref{LA.7} we obtain
  \begin{align*}
   & I \leq  \|\varphi_{u_1}\|_{L^4(\Gamma)} \| z_{u_1,v}-z_{u_2,v} \|_{L^4(\Gamma)} \|v\|_{L^2(\Gamma)}\\
    &\leq  |\Gamma|^{1/4} \|\varphi_{u_1}\|_{L^\infty(\Gamma)}
    C_\Gamma \| z_{u_1,v}-z_{u_2,v} \|_{H^1(\Omega)} \|v\|_{L^2(\Gamma)}\\
    &\leq  |\Gamma|^{1/4}C_\Gamma K^*_\infty L_{G',u}\|u_1-u_2\|_{L^s(\Gamma)} \|v\|_{L^2(\Gamma)}^2
    \leq |\Gamma|^{1/4}C_\Gamma K^*_\infty L_{G',u}2\rho_1\|v\|_{L^2(\Gamma)}^2
    = \frac{\varepsilon}{4}\|v\|_{L^2(\Gamma)}^2.
  \end{align*}
  Set $\rho_2 = \frac{\varepsilon}{8C_G C_\Gamma^2 L_{\Phi,u}}$. If $u_1,u_2\in B_{\rho_2}(u)$, using Cauchy inequality, \eqref{E05}, Lemma \ref{LA.5}, and Lemma \ref{LA.8} we get
  \begin{align*}
    II&\leq \|\varphi_{u_1}-\varphi_{u_2}\|_{L^4(\Gamma)}
    \|z_{u_2,v} \|_{L^4(\Gamma)} \|v\|_{L^2(\Gamma)}
    \leq  L_{\Phi,u} C_\Gamma \|u_1-u_2\|_{L^s(\Gamma)} C_\Gamma C_G \|v\|_{L^2(\Gamma)}^2\\
    &\leq  L_{\Phi,u} C_\Gamma^2 2\rho_2 C_G \|v\|_{L^2(\Gamma)}^2  = \frac{\varepsilon}{4}\|v\|_{L^2(\Gamma)}^2.
  \end{align*}
  Given $\varepsilon_3=\frac{\varepsilon}{8K_\infty|\Gamma|^{1/4}C_\Gamma}$, we infer from Lemma \ref{LA.10} that there exists $\rho_3>0$ such that $\|\eta_{u_1,v}-\eta_{u_2,v}\|_{H^1(\Omega)} \leq \varepsilon_3 \|v\|_{L^2(\Gamma)}$ $\forall u_1,u_2\in B_{\rho_3}(u)$. This leads to
  \begin{align*}
III   \leq &  |\Gamma|^{1/4} \|y_{u_1}\|_{L^\infty(\Omega)} \| \eta_{u_1,v} - \eta_{u_2,v}\|_{L^4(\Gamma)} \|v\|_{L^2(\Gamma)}\\
\leq  & |\Gamma|^{1/4} K_\infty C_\Gamma \varepsilon_3 \|v\|_{L^2(\Gamma)}^2 =  \frac{\varepsilon}{4}\|v\|_{L^2(\Gamma)}^2.
  \end{align*}

To estimate $IV$ we take $\rho_4 = \frac{\varepsilon}{8L_{G,u}|\Gamma|^{1/4}C_\Gamma C_\Phi}$. Then, we have with Lemmas \ref{LA.3} and \ref{LA.9} that for all $u_1,u_2\in B_{\rho_4}(u)$
  \begin{align*}
 IV \leq & |\Gamma|^{1/4}\|(y_{u_1}-y_{u_2}\|_{L^\infty(\Omega)}
    \|\eta_{u_2,v}\|_{L^4(\Gamma)} \|v\|_{L^2(\Gamma)}\\
\leq  & |\Gamma|^{1/4}L_{G,u} \|u_1-u_2\|_{L^s(\Gamma)} C_\Phi C_\Gamma \| v\|_{L^2(\Gamma)}^2 <  \frac{\varepsilon}{4}\|v\|_{L^2(\Gamma)}^2.
  \end{align*}
  And the proof concludes taking $\rho = \min\{\rho_i,\ i=0,\ldots,4\}$.
\end{proof}

\end{document}